# Domains of attraction of the random vector $(X, X^2)$ and applications

By Edward Omey

EHSAL, Stormstraat 2, B-1000 Brussels - Belgium

edward.omey@ehsal.be; www.edwardomey.com

## Abstract

Many statistics are based on functions of sample moments. Important examples are the sample variance $s_{n-1}^2$, the sample coefficient of variation $SV(n)$, the sample dispersion $SD(n)$ and the non-central $t$-statistic $t(n)$. The definition of these quantities makes clear that the vector defined by

$$(\sum_{i=1}^{n} X_i, \sum_{i=1}^{n} X_i^2)$$

plays an important role. In studying the asymptotic behaviour of this vector we start by formulating best possible conditions under which the vector $(X, X^2)$ belongs to a bivariate domain of attraction of a stable law. This approach is new, uniform and simple. Our main results include a full discussion of the asymptotic behaviour of $SV(n)$, $SD(n)$ and $t^2(n)$. For simplicity, in restrict ourselves to positive random variables $X$.

**Keywords and Phrases**

Domains of attraction, regular variation, sample dispersion, sample coefficient of variation, non-central $t$-statistic

# 1  Introduction

Let $F(x) = P(X \leq x)$ denote the distribution function (d.f.) of a positive random variable $X$ and let $\overline{F}(x) = 1 - F(x)$ denote the tail. Let $G(x,y)$ denote the d.f. of the random vector (r.v.) $(X, X^2)$. Using this notation we find that $G(x,y) = P(X \leq x, X^2 \leq y) = F(\min(x, \sqrt{y}))$. This relationship can be exploited to transfer many properties from $F$ to $G$. Studying the random vector $(X, X^2)$ can be interesting because it is linked to many statistical estimators. To study the mean $\mu = E(X)$ and the variance $\sigma^2 = Var(X)$ for example, one uses the sample mean $\overline{X}$ and the sample variance $s_{n-1}^2 = n(\overline{X^2} - \overline{X}^2)/(n-1)$. Other related statistical measures are the non-central $t$-statistic $t(n) = \sqrt{n}\,\overline{X}/s$, the coefficient of variation $SV(n) = s/\overline{X}$ and the sample dispersion $SD(n) = s^2/\overline{X}$. Many asymptotic properties of these statistics are known if the mean and the variance are finite. On the other hand, if the variance or the mean is not finite, it also makes sense to study asymptotic properties of these quantities. In the case of the statistic $t(n)$ and $\mu = 0$, the statistic $t(n)$ is called the central $t$-statistics and its limiting behaviour is now well understood. See, for example [18], [4], [7]. If $\mu \neq 0$, the statistic $t(n)$ is called the non-central $t$-statistic. In contrast to the central $t$-statistic, there are only a few studies on the limiting behaviour of $t(n)$. For a recent paper devoted to $t(n)$, we refer to [2] and the extensive list of references given there. In [1] and ]16], the authors discuss asymptotic properties of $SV(n)$ and $SD(n)$ and. They also consider $SV(N(t))$ and $SD(N(t))$ where $N(t)$ is an integer valued random variable. It turns out, cf. section 4, that a key role is played by the vector $(\sum_{i=1}^{n} X_i, \sum_{i=1}^{n} X_i^2)$ and by $T(n)$, where

$$T(n) = \frac{\sum_{i=1}^{n} X_i^2}{(\sum_{i=1}^{n} X_i)^2}$$

In the present paper we offer a unified approach towards the statistics mentioned above. As the starting point of our analysis we discuss conditions under which the vector $(X, X^2)$ belongs to the bivariate domain of attraction of a bivariate stable law. By applying simple transformations, we then obtain the desired results.

The paper is organised as follows. In section 2 we briefly discuss univariate and bivariate domains of attraction. In section 3 we discuss domains of attraction of the vector $(X, X^2)$ and in section 4 we use transformations to recover many results concerning $SV(n)$, $SD(n)$ and $t(n)$. We finish the paper with some concluding remarks.

In the paper we restrict ourself to nonnegative random variables $X$. The case of real $X$ will be postponed to another paper. In the real case one can prove similar results, but one has to study different cases depending on whether $\mu = 0$ or $\mu \neq 0$.

For further use, recall the definition of regular variation. A positive and measurable function $g(x)$ is regularly varying with real index $\alpha$ (notation $g \in RV(\alpha)$) if as $t \to \infty$, $g(ty)/g(t) \to y^\alpha$, $\forall y > 0$. It can be proved that the defining convergence holds locally uniformly (l.u.) with respect to $y > 0$. For this and other properties and applications of regular variation, we refer to [3], [6] or [22]. For a recent survey paper, see [12]. For applications in extreme value theory we refer to [11]. A point process approach to regular variation and weak convergence is available in [20].

## 2 Domains of attraction

In this section we briefly discuss univariate and bivariate domains of attraction of nonnegative random variables and vectors. These results are not new and are well known by now.

### 2.1 Univariate case

Recall that the random variable $X$ belongs to the domain of attraction of a stable law $Y(\alpha)$ with parameter $\alpha$, $0 < \alpha \leq 2$, if there exist positive numbers $a(n)$ and real numbers $c(n)$ so that

$$\lim_{n \to \infty} P(\frac{S(n) - c(n)}{a(n)} \leq x) = P(Y(\alpha) \leq x)$$

Notation $X \in D(Y(\alpha))$. Here $S(n) = X_1 + X_2 + ... + X_n$ is the sequence of partial sums generated by i.i.d. copies of $X$. Note that for $\alpha = 2$, $Y(2)$ has a normal distribution. Now assume that $X \geq 0$ and let $F(x) = P(X \leq x)$. For further use we define the truncated second moment function $V(x)$:

$$V(x) = \int_0^x y^2 dF(y)$$

The following result is well known, cf. [5], [19].

**Proposition 1** *Suppose that $X$ is nonnegative and not concentrated in 1 point.*

*(i) For $0 < \alpha < 2$, we have $X \in D(Y(\alpha))$ if and only if $\overline{F} \in RV(-\alpha)$.*

*(ii) For $\alpha = 2$ we have $X \in D(Y(2))$ if and only if $V \in RV(0)$.*

Note that for $0 < \alpha \leq 2$, $X \in D(Y(\alpha))$ is equivalent to
$$x^2 \overline{F}(x)/V(x) \to (2-\alpha)/\alpha.$$

Moreover, $\overline{F} \in RV(-2)$ implies that $V \in RV(0)$ and then $X \in D(Y(2))$.

We have some freedom in choosing the normalizing sequences $\{a(n)\}$ and $\{c(n)\}$. In our paper, we use the normalizing constants by replacing $x$ by $n$ in the functions $a(x)$ and $c(x)$ defined as follows:

- For $\alpha = 2$, $a(x) \in RV(1/2)$ is determined by the asymptotic relation
$$xV(a(x))/a^2(x) \to 1. \tag{1}$$

- If $0 < \alpha < 2$, $a(x) \in RV(1/\alpha)$ is determined by by the asymptotic relation
$$x\overline{F}(a(x)) \to 1. \tag{2}$$

- If $0 < \alpha < 1$, we choose $c(x) = 0$. If $1 < \alpha \leq 2$, we choose $c(x) = x\mu = xE(X)$.

- If $\alpha = 1$, then $c(x) \in RV(1)$ is given by the relation $c(x) = xm(a(x))$ where $a(x)$ is given by (2) and where $m(x)$ denotes the integrated tail
$$m(x) = \int_0^x \overline{F}(t)dt.$$

Note that if $\mu_2 = EX^2 < \infty$, then $V(x) \to \mu_2$ and $a^2(n) \sim n\mu_2$. To obtain the precise form of the limit $Y(\alpha)$, we denote by $\widehat{g}(s)$ the Laplace-Stieltjes transform of a function $G(x)$:
$$\widehat{G}(s) = \int_0^\infty \exp(-sx)dG(x).$$

If $G(x) = P(Z \leq x)$ is the distribution function of $Z \geq 0$, we have $\widehat{G}(s) = E(\exp -sZ)$, the generating function of $Z$. For $F(x)$ we write $\varphi(s) = \widehat{F}(s) = E(\exp(-sX))$ and in the results below we use the notation $\varphi_\alpha(s) = E(\exp -sY(\alpha))$. With our choice of $\{a(n)\}$ and $\{c(n)\}$ we have:

(i) If $0 < \alpha < 1$ we have $\phi_\alpha(s) = \exp(-\Gamma(1-\alpha)s^\alpha)$;
(ii) If $\alpha = 1$ we have $\phi_\alpha(s) = \exp(-s\log(s) - \gamma s)$;
(iii) If $1 < \alpha < 2$, we have $\phi_\alpha(s) = \exp(-\Gamma(2-\alpha)s^\alpha/(\alpha-1))$;
(iv). If $\alpha = 2$, we have $\phi_\alpha(s) = \exp(\frac{1}{2}s^2(1-\mu^2/V(\infty)))$.

Note that if $V(\infty) = \infty$ in (iv) we have $\phi_\alpha(s) = \exp(\frac{1}{2}s^2)$. If $\mu_2 < \infty$, then $\phi_\alpha(s) = \exp(\frac{1}{2}s^2 Var(X)/\mu_2)$.

**Remark.** If $1 \leq \alpha < 2$, then we have $S(n)/nm(a(n)) \xrightarrow{P} 1$.

## 2.2 Multivariate case

Many of the one-dimensional results can and have been generalized to higher dimensions. Recall that the random vector $(X, Y)$ belongs to the bivariate domain of attraction of a stable vector $(Y_1(\alpha), Y_2(\beta))$ if we can find sequences of constants $a(n) > 0, b(n) > 0$ and $c(n), d(n)$ such that

$$P(\frac{S_X(n) - c(n)}{a(n)} \leq x, \frac{S_Y(n) - d(n)}{b(n)} \leq y) \to P(Y_1(\alpha) \leq x, Y_2(\beta) \leq y) \quad (3)$$

where $S_X(n)$ and $S_Y(n)$ are the partial sums of independent copies of $(X, Y)$. Notation $(X, Y) \in D(Y_1(\alpha), Y_2(\beta))$. Assuming that $Y_1(\alpha)$ and $Y_2(\beta)$ are nondegenerate, the normalizing constants are determined by the convergence of the marginals in (3) and then we use the normalizing constants given as in Section 2.1. We denote the d.f. of $(X, Y)$ by $F(x, y)$. For $(Y_1(\alpha), Y_2(\beta))$, we denote the generating function by $\psi_{1,2}(s, t)$ and its marginals $\psi_1(s)$, $\psi_2(t)$.

We assume that $(X, Y) \geq (0, 0)$ and define the integrated "tail" $U(x, y)$ function as:

$$U(x, y) = \int_0^x \int_0^y P(X > u, Y > v) du dv.$$

A bivariate analogue of the truncated "second moment" $V(x)$ can be defined as $W(x, y)$, where

$$W(x, y) = \int_0^x \int_0^y uv dF(u, v).$$

If $E(XY) < \infty$, we see that $W(x, y) \to E(XY)$ as $\min(x, y) \to \infty$. Using partial integration, it easily follows that

$$|W(x, y) - U(x, y)| \leq 2xy P(X > x) + 2xy P(Y > y) \quad (4)$$

and (4) can be used to transfer properties of $U$ to $W$ and vice versa. In applications one can choose the more convenient function.

The next result goes back to [21]. We state a simplified version of a result of [8], see also [9], [10].

**Theorem 2** *Suppose $(X, Y) \geq (0, 0)$ and $X \in D(Y_1(\alpha))$, $Y \in D(Y_2(\beta))$ with $0 < \alpha, \beta \leq 2$.*

*(i) If $0 < \alpha, \beta < 2$, then $(X, Y) \in D(Y_1(\alpha), Y_2(\beta))$ iff for each $x, y > 0$ we have*

$$\lim_{t \to \infty} \frac{tU(a(t)x, b(t)y)}{a(t)b(t)} = \lim_{y \to \infty} \frac{tW(a(t)x, b(t)y)}{a(t)b(t)} = \Omega(x, y) \quad (5)$$

*for some limit function* $\Omega(x,y)$. *In this case we have*

$$\psi_{1,2}(s,t) = \psi_1(s)\psi_2(t)\exp st\widehat{\Omega}(s,t).$$

(ii) *If* $\alpha = \beta = 2$, *then* $(X,Y) \in D(Y_1(2), Y_2(2))$ *iff for each* $x,y > 0$ *we have*

$$\lim_{t\to\infty} \frac{tU(a(t)x, b(t)y)}{a(t)b(t)} = \lim_{t\to\infty} \frac{tW(a(t)x, b(t)y)}{a(t)b(t)} = C,$$

*where* $C \geq 0$ *is a constant. In this case we have*

$$\psi_{1,2}(s,t) = \psi_1(s)\psi_2(t)\exp(C-c)st,$$

*where* $c = E(X)E(Y)/\sqrt{E(X^2)E(Y^2)}$ *if* $E(X^2 + Y^2) < \infty$ *and* $c = 0$ *otherwize.*

(iii) *If* $0 < \alpha < 2$ *and* $\beta = 2$, *then* $(X,Y) \in D(Y_1(\alpha), Y_2(2))$ *and* $Y_1(\alpha)$ *and* $Y_2(2)$ *are independent.*

**Remarks 2.**

1) In the case where $0 < \alpha, \beta < 2$, one can prove (c. [10]) that (5) can be replaced by

$$\lim_{t\to\infty} t\overline{F}_{X,Y}(a(t)x, b(t)y) = h(x,y), \quad \forall x,y > 0,$$

for some finite limit function $h(x,y)$.

2) In Theorem 2(ii) the limiting bivariate normal distribution has variance-covariance matrix given by

$$\Sigma = \begin{bmatrix} 1 - E^2(X)/E(X^2) & C - c \\ C - c & 1 - E^2(Y)/E(Y^2) \end{bmatrix}$$

where $E^2(X)/E(X^2) = 0$ if $E(X^2) = \infty$. and $E^2(Y)/E(Y^2) = 0$ if $E(Y^2) = \infty$. If $E(X^2 + Y^2) = \infty$, we have $c = 0$. If $E(X^2 + Y^2) < \infty$, we have $E(XY) < \infty$ and then we have

$$C = \frac{E(XY)}{\sqrt{E(X^2)E(Y^2)}}.$$

In this case the limiting normal distribution has variance-covariance matrix given by

$$\Sigma = \begin{bmatrix} Var(X)/E(X^2) & Cov(X,Y)/\sqrt{E(X^2)E(Y^2)} \\ Cov(X,Y)/\sqrt{E(X^2)E(Y^2)} & Var(Y)/E(Y^2) \end{bmatrix}$$

# 3  The case where $\mathbf{X} = (X, X^2)$

For the vector $(X, X^2)$ with $X \geq 0$, we have $F(x) = P(X \leq x)$ and $F_2(x) = P(X^2 \leq x) = F(\sqrt{x})$. Clearly $\overline{F} \in RV(-\alpha)$ holds if and only if $\overline{F}_2 \in RV(-\alpha/2)$. In the next result we use the notation $\mu = E(X)$ and $\mu_k = E(X^k)$ whenever needed.

**Theorem 3** *(i) For $0 < \alpha < 2$ we have $(X, X^2) \in D(Y_1(\alpha), Y_2(\alpha/2))$ iff $\overline{F} \in RV(-\alpha)$.*
  *(ii) For $2 \leq \alpha < 4$ we have $(X, X^2) \in (Y_1(2), Y_2(\alpha/2))$ iff $\overline{F} \in RV(-\alpha)$. Moreover, $Y_1(2)$ and $Y_2(\alpha/2)$ are independent.*
  *(iii) We have $(X, X^2) \in (Y_1(2), Y_2(2))$ iff $X^2 \in D(Y_2(2))$.*

**Proof.** (i) For $0 < \alpha < 2$, the univariate results show that $\overline{F} \in RV(-\alpha)$ if and only if $X \in D(Y_1(\alpha))$ and $X^2 \in D(Y_2(\alpha/2))$. Using (2) we see that $a^2(t) = b(t)$. Now we consider $U(x, y)$. Clearly we have

$$U(x, y) = \int_0^x \int_0^y \overline{F}_X(\max(u, \sqrt{v})) du dv.$$

and

$$\frac{U(tx, t^2 y)}{t^3 \overline{F}_X(t)} = \int_0^x \int_0^y \frac{\overline{F}_X(t \max(u, \sqrt{v}))}{\overline{F}_X(t)} du dv.$$

Using $\overline{F}_X(t \max(u, \sqrt{v})) \leq \overline{F}_X(t\sqrt{v})$ and $\overline{F}(\sqrt{x}) \in RV(-\alpha/2)$ we can apply known properties of regular variation to obtain that

$$\frac{U(tx, t^2 y)}{t^3 \overline{F}_X(t)} \to \Omega(x, y) = \int_0^x \int_0^y (\max(u, \sqrt{v}))^{-\alpha} du dv.$$

Replacing $t$ by $a(t)$ we find that $a^3(t) \overline{F}(a(t)) \sim a^3(t)/t$ and Theorem 2 shows that $(X, X^2) \in D((Y_1(\alpha), Y_2(\alpha/2))$.
  (ii) Note that for $2 \leq \alpha < 4$, $\overline{F} \in RV(-\alpha)$ implies that $X \in D(Y_1(2))$. Also $\overline{F} \in RV(-\alpha)$ holds if and only if $X^2 \in D(Y_2(\alpha/2))$. Now Theorem 2 can be used.
  (iii) If $X^2 \in D(Y_2(2))$ then automatically we have $X \in D(Y_1(2))$ and $\mu_3 = W(\infty, \infty) < \infty$. It follows that

$$\lim \frac{tW(a(t)x, a^2(t)y)}{a(t)a^2(t)} = C$$

where $C = 0$ if $\mu_4 = \infty$ and $C = \mu_3/\sqrt{\mu_2 \mu_4}$ if $\mu_4 < \infty$. Theorem 2 applies and we find that $(X, X^2) \in D(Y_1(2), Y_2(2))$. ∎

**Remarks 3.**.

1) In case (iii) of Theorem 3, the limiting normal law has the variance-covariance matrix given by

$$\Sigma = \begin{bmatrix} (\mu_2 - \mu^2)/\mu_2 & 0 \\ 0 & 1 \end{bmatrix}$$

if $\mu_4 = \infty$ and by

$$\Sigma = \begin{bmatrix} (\mu_2 - \mu^2)/\mu_2 & (\mu_3 - \mu_1\mu_2)/\sqrt{\mu_2\mu_4} \\ (\mu_3 - \mu_1\mu_2)/\sqrt{\mu_2\mu_4} & (\mu_4 - \mu_2^2)/\mu_4 \end{bmatrix}$$

if $\mu_4 < \infty$.

2) Theorem 3 gives conditions under which

$$(\frac{\sum_{i=1}^n X_i - c(n)}{a(n)}, \frac{\sum_{i=1}^n X_i^2 - d(n)}{b(n)}) \overset{d}{\Longrightarrow} (Y_1(u), Y_2(v)) \qquad (6)$$

for some numbers $u$ and $v$. For further use, in this remark we give the precise form of the normalizing sequences.

(i) If $0 < \alpha < 2$, then in (6) we have $u = \alpha$ and $v = \alpha/2$. We can use (2) to see that $b(n) = a^2(n)$. Since $\alpha/2 < 1$, we can take $d(n) = 0$ and $c(n)$ according to the different cases of Section 2.1.

(ii) If $2 \leq \alpha < 4$, then in (6) we have $u = 2$ and $v = \alpha/2$. We have $c(n) = n\mu$ and $a(n)$ is determined by (1). The sequence $b(n)$ is determined by the relation $n(1 - F_2(b(n))) \to 1$. If $\alpha > 2$ we have $d(n) = n\mu_2$ while if $\alpha = 2$ we take

$$d(n) = nm_2(b(n)) \qquad (7)$$

where $m_2(x) = \int_0^x \overline{F}_2(u)du$.

(iii) In case (iii) of Theorem 3, in (6) we have $u = v = 2$. Now we take $c(n) = n\mu$ and $d(n) = n\mu_2$. The sequence $a(n)$ is determined by (1) and $b(n)$ is determined by

$$\frac{nV_2(b(n))}{b^2(n)} \to 1 \qquad (8)$$

where $V_2(x) = E(X^4 I\{X^2 \leq x\}) = \int_0^{\sqrt{x}} u^4 dF(u)$.

## 4 Applications

As mentioned in the introduction, many characteristics in statistics are based on $(\sum_{i=1}^n X_i, \sum_{i=1}^n X_i^2)$. As examples we mention the sample coef-ficient of variation $SV(n)$, the sample dispersion $SD(n)$ and $t(n)$. As in [16], we define $C(n)$ and $T(n)$ as follows:

$$C(n) = \frac{\sum_{i=1}^n X_i^2}{\sum_{i=1}^n X_i}, \; T(n) = \frac{\sum_{i=1}^n X_i^2}{(\sum_{i=1}^n X_i)^2}.$$

Note that
$$SV(n) = \sqrt{nT(n) - 1}$$
$$SD(n) = C(n) - \overline{X} = (nT(n) - 1)\overline{X}$$

and
$$t(n) = \frac{\sqrt{n}}{\sqrt{nT(n) - 1}}.$$

Clearly $T(n)$ plays an important role in studying $SV(n)$ and $t(n)$.

## 4.1 Asymptotic behaviour of $T(n)$ and $SV(n)$

Using the notations of Theorem 2 and the Remark 3.2. we have the following result. The result simplifies and completes the proofs of the corresponding results of [1], [16], [17]. Also cases (iii) and (v) seem to be new.

**Theorem 4** *(i) Suppose that $\overline{F}(x) \in RV(-\alpha)$ with $0 < \alpha < 1$. Then*

$$T(n) \overset{d}{\Longrightarrow} \frac{Y_2(\alpha/2)}{Y_1^2(\alpha)}.$$

*(ii) Suppose that $\overline{F}(x) \in RV(-\alpha)$ with $1 \leq \alpha < 2$. Then*

$$\frac{n^2 m^2(a(n))}{a^2(n)} T(n) \overset{d}{\Longrightarrow} Y_2(\alpha/2).$$

*(iii) Suppose that $\overline{F}(x) \in RV(-2)$. Then*

$$\frac{n}{b(n)}(nT(n) - \frac{d(n)}{n\mu^2}) \overset{d}{\Longrightarrow} \frac{1}{\mu^2} Y_2(1)$$

*where $d(n)$ is given in (7).*
*(iv) Suppose that $\overline{F}(x) \in RV(-\alpha)$ with $2 < \alpha < 4$. Then*

$$\frac{n}{b(n)}(nT(n) - \frac{\mu_2}{\mu^2}) \overset{d}{\Longrightarrow} \frac{1}{\mu^2} Y_2(\alpha/2).$$

(v) Suppose that $X^2 \in D(Y_2(2))$. Then
$$\frac{n}{b(n)}(nT(n) - \frac{\mu_2}{\mu^2}) \stackrel{d}{\Longrightarrow} Y_3(2)$$

where
$$Y_3(2) \stackrel{d}{=} \frac{1}{\mu^2}Y_2(2) - 2\frac{\mu_2\sqrt{\mu_2}}{\mu^3\sqrt{\mu_4}}Y_1(2)$$

if $\mu_4 < \infty$ and $Y_3(2) \stackrel{d}{=} Y_2(2)/\mu^2$ otherwize.

**Proof.** (i) Under the conditions of (i) we have
$$(\frac{\sum_{i=1}^n X_i}{a(n)}, \frac{\sum_{i=1}^n X_i^2}{a^2(n)}) \stackrel{d}{\Longrightarrow} (Y_1(\alpha), Y_2(\alpha/2)).$$

Now it follows that
$$P(T(n) \leq x) = P(\frac{1}{a^2(n)}\sum_{i=1}^n X_i^2 \leq x(\frac{1}{a(n)}\sum_{i=1}^n X_i)^2)$$

so that
$$P(T(n) \leq x) \to P(Y_2(\alpha/2) \leq xY_1^2(\alpha)).$$

(ii) Now we have
$$(\frac{\sum_{i=1}^n X_i}{nm(a(n))}, \frac{\sum_{i=1}^n X_i^2}{a^2(n)}) \stackrel{d}{\Longrightarrow} (1, Y_2(\alpha/2))$$

and the result follows as in (i).

(iii) Using the notations
$$A(n) = \frac{(\sum_{i=1}^n X_i)^2 - (n\mu)^2}{na(n)}, \quad B(n) = \frac{\sum_{i=1}^n X_i^2 - d(n)}{b(n)}.$$

in this case (iii) we have
$$(A(n), B(n)) \stackrel{d}{\Longrightarrow} (2\mu Y_1(2), Y_2(1)). \tag{9}$$

To prove the result, we consider
$$I = P(\frac{d(n)}{b(n)}(\frac{n^2 T(n)}{d(n)} - \frac{1}{\mu^2}) \leq x).$$

Using the definition of $T(n)$, we obtain that

$$\begin{aligned} I &= P(\sum_{i=1}^{n} X_i^2 - (\frac{d(n)}{\mu^2 n^2} + \frac{xb(n)}{n^2})(\sum_{i=1}^{n} X_i)^2 \leq 0) \\ &= P(B(n)b(n) - (\frac{d(n)}{\mu^2 n^2} + \frac{xb(n)}{n^2})A(n)na(n) \leq xb(n)\mu^2) \\ &= P(B(n) - (\frac{a(n)d(n)}{\mu^2 b(n)n} + \frac{a(n)x}{n})A(n) \leq x\mu^2). \end{aligned}$$

Now recall that $A(n) \overset{d}{\Longrightarrow} 2\mu Y_2(2)$ and observe that $a(x) \in RV(1/2)$. Also note that $b(x) \in RV(1)$ and (cf. (7)) that $d(x) \in RV(1)$. It follows that $a(x)/x \to 0$ and that $a(x)d(x)/(xb(x)) \to 0$. Using (9), we conclude that $I \to P(Y_2(1) \leq x\mu^2)$.

(iv) In this case we have $d(n) = n\mu_2$ and

$$(A(n), B(n)) \overset{d}{\Longrightarrow} (2\mu Y_1(2), Y_2(\alpha/2)) \tag{10}$$

where, cf. (1), $a^2(n) \sim n\mu_2$ and $b(x) \in RV(2/\alpha)$. Now consider

$$I = P(\frac{n}{b(n)}(nT(n) - \frac{\mu_2}{\mu^2}) \leq x).$$

By using the definition of $T(n)$, we find that

$$\begin{aligned} I &= P(\sum_{i=1}^{n} X_i^2 - (\sum_{i=1}^{n} X_i)^2(\frac{xb(n)}{n^2} + \frac{\mu_2}{n\mu^2}) \leq 0) \\ &= P(B(n)b(n) - (\frac{xb(n)}{n^2} + \frac{\mu_2}{n\mu^2})na(n)A(n) \leq xb(n)\mu^2) \\ &= P(B(n) - (\frac{a(n)x}{n} + \frac{a(n)\mu_2}{b(n)\mu^2})A(n) \leq \mu^2 x). \end{aligned}$$

Since $a(n)/n \to 0$ and $a(n)/b(n) \to 0$ as $n \to \infty$, Using (10), we conclude that $I \to P(Y_2(\alpha/2) - 0 \leq \mu^2 x)$ and this proves the result.

(v) Now we have

$$(A(n), B(n)) \overset{d}{\Longrightarrow} (2\mu Y_1(2), Y_2(2))$$

where, cf. (1), $a^2(n) \sim n\mu_2$ and $b(n)$ is determined by (8). As in case (iv) we consider $I$ and again we find that

$$I = P(B(n) - (\frac{a(n)x}{n} + \frac{a(n)\mu_2}{b(n)\mu^2})A(n) \leq \mu^2 x)$$

and now we have to distinghuish between two cases. If $\mu_4 < \infty$, we find that $b^2(n) \sim n\mu_4$ and it follows that

$$I \to P(Y_2(2) - \frac{\mu_2\sqrt{\mu_2}}{\mu^2\sqrt{\mu_4}} 2\mu Y_1(2) \leq \mu^2 x).$$

In the case where $\mu_4 = \infty$, we have

$$\frac{a^2(n)}{b^2(n)} \sim \frac{n\mu_2}{b^2(n)} \sim \frac{\mu_2}{V_2(b(n))} \to 0.$$

and we find that $I \to P(Y_2(2) \leq \mu^2 x)$. ∎

**Remarks 4.**

1) Using the variance-covariance matrix $\Sigma$ one can calculate the variance of the limiting normal random variable $Y_3(2)$ in Theorem 4(v).

2) The prominent place of $\mu$ in the theorem shows that we can expect difficulties if $X$ is a real random variable with $\mu = 0$.

3) In the special case where $P(X = 1) = p$, $P(X = 0) = q$, with $0 < p = 1 - q < 1$, we find that

$$\frac{T(n)}{n} = \frac{1}{\overline{X}}$$

and we obtain that

$$\sqrt{n}(\frac{T(n)}{n} - \frac{1}{p}) = \frac{\sqrt{n}(p - \overline{X})}{p\overline{X}} \xrightarrow{d} -\frac{\sqrt{pq}}{p^2} Z$$

where $Z \sim N(0,1)$. Theorem 4(v) (where we used different normalizing sequences) gives the same result taking into account that now we have $Y_1(2) = Y_2(2)$.

Now we can use Theorem 4 to obtain the precise asymptotic behaviour of the sample coefficient of variation $SV(n) = \sqrt{nT(n) - 1}$. The next result completes the results of [16, Section 4.1]. We use the notation $\sigma^2 = \mu_2 - \mu^2$.

**Corollary 5** *(i) Suppose that $\overline{F}(x) \in RV(-\alpha)$ with $0 < \alpha < 1$. Then*

$$\frac{SV(n)}{\sqrt{n}} \xrightarrow{d} \frac{\sqrt{Y_2(\alpha/2)}}{Y_1(\alpha)}.$$

*(ii) Suppose that $\overline{F}(x) \in RV(-\alpha)$ with $1 \leq \alpha < 2$. Then*

$$\frac{\sqrt{n}m(a(n))}{a(n)} SV(n) \xrightarrow{d} \sqrt{Y_2(\alpha/2)}.$$

(iii) Suppose that $\overline{F}(x) \in RV(-2)$. Then

$$\frac{n\sqrt{c(n)}}{b(n)}(SV(n) - \sqrt{c(n)}) \stackrel{d}{\Longrightarrow} \frac{1}{2\mu^2}Y_2(1)$$

where $c(n)$ is given by $c(n) = d(n)/n\mu^2 - 1$.

(iv) Suppose that $\overline{F}(x) \in RV(-\alpha)$ with $2 < \alpha < 4$. Then

$$\frac{n}{b(n)}(SV(n) - \frac{\sigma}{\mu}) \stackrel{d}{\Longrightarrow} \frac{1}{2\sigma\mu}Y_2(\alpha/2).$$

(v) Suppose that $X^2 \in D(Y_2(2))$. Then

$$\frac{n}{b(n)}(SV(n) - \frac{\sigma}{\mu}) \stackrel{d}{\Longrightarrow} \frac{\mu}{2\sigma}Y_3(2)$$

where $Y_3(2)$ is given in Theorem 4(v).

**Proof.** (i) and (ii) follow immediately from Theorem 4 (i),(ii).

(iii) From Theorem 4(iii) we obtain that

$$\frac{nT(n)}{d(n)/n} \stackrel{p}{\longrightarrow} \frac{1}{\mu^2}$$

and

$$\frac{n}{b(n)}(nT(n) - 1 - c(n)) \stackrel{d}{\Longrightarrow} \frac{1}{\mu^2}Y_2(1)$$

where $c(n)$ is given as in (iii). It follows that $SV(n)/\sqrt{c(n)} \stackrel{p}{\longrightarrow} 1$. Now observe that

$$SV(n) - \sqrt{c(n)} = \frac{nT(n) - 1 - c(n)}{SV(n) + \sqrt{c(n)}}.$$

Using Theorem 4(iii), we obtain that

$$\frac{n\sqrt{c(n)}}{b(n)}(SV(n) - \sqrt{c(n)}) = \frac{1}{2\mu^2}Y_2(1).$$

(iv) First note that Theorem 4(iv) shows that $nT(n) \stackrel{p}{\longrightarrow} \mu_2/\mu^2$ so that $SV(n) \stackrel{p}{\longrightarrow} \sigma/\mu$ where $\sigma^2 = \mu_2 - \mu^2$. Now we have

$$SV(n) - \frac{\sigma}{\mu} = \frac{nT(n) - \mu_2/\mu^2}{\sigma/\mu + SV(n)}.$$

Theorem 4(iv) can now be used to obtain the desired result.

(v) The proof is similar and therefore omitted. ∎

## 4.2 Asymptotic behaviour of $C(n)$ and $SD(n)$

Following similar steps as in the proof of Theorem 4 we readily obtain the following results for $C(n)$.

**Theorem 6** *(i) Suppose that $\overline{F}(x) \in RV(-\alpha)$ with $0 < \alpha < 1$. Then*

$$\frac{1}{a(n)} C(n) \stackrel{d}{\Longrightarrow} \frac{Y_2(\alpha/2)}{Y_1(\alpha)}.$$

*(ii) Suppose that $\overline{F}(x) \in RV(-\alpha)$ with $1 \leq \alpha < 2$. Then*

$$\frac{nm(a(n))}{a^2(n)} C(n) \stackrel{d}{\Longrightarrow} Y_2(\alpha/2).$$

*(iii) Suppose that $\overline{F}(x) \in RV(-2)$. Then*

$$\frac{d(n)}{b(n)} \left( \frac{n}{d(n)} C(n) - \frac{1}{\mu} \right) \stackrel{d}{\Longrightarrow} \frac{1}{\mu} Y_2(1).$$

*(iv) Suppose that $\overline{F}(x) \in RV(-\alpha)$ with $2 < \alpha < 4$. Then*

$$\frac{n}{b(n)} \left( C(n) - \frac{\mu_2}{\mu} \right) \stackrel{d}{\Longrightarrow} \frac{1}{\mu} Y_2(\alpha/2).$$

*(v) Suppose that $X^2 \in D(Y_2(2))$. Then*

$$\frac{n}{b(n)} \left( C(n) - \frac{\mu_2}{\mu} \right) \stackrel{d}{\Longrightarrow} Y_3(2)$$

*where*

$$Y_3(2) \stackrel{d}{=} \frac{1}{\mu} Y_2(2) - \frac{\mu_2 \sqrt{\mu_2}}{\mu^2 \sqrt{\mu_4}} Y_1(2)$$

*if $\mu_4 < \infty$ and $Y_3(2) \stackrel{d}{=} Y_2(2)/\mu$ otherwize.*

Now we can use Theorem 6 to discuss the asymptotic behaviour of $SD(n) = C(n) - \overline{X}$. The proof of the next result is left for the reader.

**Corollary 7** *(i) Suppose that $\overline{F}(x) \in RV(-\alpha)$ with $0 < \alpha < 1$. Then*

$$\frac{1}{a(n)} SD(n) \stackrel{d}{\Longrightarrow} \frac{Y_2(\alpha/2)}{Y_1(\alpha)}.$$

(ii) Suppose that $\overline{F}(x) \in RV(-\alpha)$ with $1 \leq \alpha < 2$. Then

$$\frac{nm(a(n))}{a^2(n)} SD(n) \overset{d}{\Longrightarrow} Y_2(\alpha/2).$$

(iii) Suppose that $\overline{F}(x) \in RV(-2)$. Then

$$\frac{n}{b(n)}(SD(n) - \frac{d(n)}{n\mu} + \mu) \overset{d}{\Longrightarrow} \frac{1}{\mu} Y_2(1).$$

(iv) Suppose that $\overline{F}(x) \in RV(-\alpha)$ with $2 < \alpha < 4$. Then

$$\frac{n}{b(n)}(SD(n) - \frac{\sigma^2}{\mu}) \overset{d}{\Longrightarrow} \frac{1}{\mu} Y_2(\alpha/2).$$

(v) Suppose that $X^2 \in D(Y_2(2))$. Then

$$\frac{n}{b(n)}(SD(n) - \frac{\sigma^2}{\mu}) \overset{d}{\Longrightarrow} Y_4(2)$$

where

$$Y_4(2) \overset{d}{=} \frac{1}{\mu} Y_2(2) - (\frac{\sigma^2}{\mu} + 2)\frac{\sqrt{\mu_2}}{\sqrt{\mu_4}} Y_1(2)$$

if $\mu_4 < \infty$ and $Y_4(2) \overset{d}{=} Y_2(2)/\mu$ otherwize.

### 4.3 The asymptotic behaviour of $t^2(n)$

Using the definition of $t(n)$ we see that $t^2(n) = n/(nT(n) - 1)$ and this relation can by used to transfer the asymptotic properties of $T(n)$ to $t^2(n)$. Our Theorem 8 covers all cases for positive random variables. It should be compared to the results of [2].

**Theorem 8** (i) Suppose that $\overline{F}(x) \in RV(-\alpha)$ with $0 < \alpha < 1$. Then

$$t^2(n) \overset{d}{\Longrightarrow} \frac{Y_1^2(\alpha)}{Y_2(\alpha/2)}.$$

(ii) Suppose that $\overline{F}(x) \in RV(-\alpha)$ with $1 \leq \alpha < 2$. Then

$$\frac{a^2(n)}{n^2 m^2(a(n))} t^2(n) \overset{d}{\Longrightarrow} \frac{1}{Y_2(\alpha/2)}.$$

(iii) Suppose that $\overline{F}(x) \in RV(-2)$. Then

$$\frac{nc^2(n)}{b(n)}(\frac{1}{c(n)} - \frac{t^2(n)}{n}) \stackrel{d}{\Longrightarrow} \frac{1}{\mu}Y_2(1),$$

where $c(n) = (d(n) - n\mu^2)/(n\mu^2)$.

(iv) Suppose that $\overline{F}(x) \in RV(-\alpha)$ with $2 < \alpha < 4$. Then

$$\frac{nc^2}{b(n)}(\frac{1}{c} - \frac{t^2(n)}{n}) \stackrel{d}{\Longrightarrow} \frac{1}{\mu^2}Y_2(\alpha/2),$$

where $c = \sigma^2/\mu^2$.

(v) Suppose that $X^2 \in D(Y_2(2))$. Then

$$\frac{nc^2}{b(n)}(\frac{1}{c} - \frac{t^2(n)}{n}) \stackrel{d}{\Longrightarrow} Y_3(2)$$

where $c = \sigma^2/\mu^2$ and $Y_3(2)$ as in Theorem 4(v).

**Proof.** (i), (ii) This follows immediately from Theorem 4(i), (ii).

(iii) To prove this result, we write

$$I = P(\frac{nc^2(n)}{b(n)}(\frac{1}{c(n)} - \frac{t^2(n)}{n}) \leq x).$$

We have

$$\begin{aligned} I &= P(\frac{t^2(n)}{n} \geq \frac{1}{c(n)} - \frac{xb(n)}{nc^2(n)}) \\ &= P(\frac{1}{nT(n)-1} \geq \frac{nc(n) - xb(n)}{nc^2(n)}) \\ &= P(\frac{n}{b(n)}(nT(n) - 1 - c(n)) \leq \frac{xnc(n)}{nc(n) - xb(n)}) \end{aligned}$$

From Theorem 4(iii) we have

$$\frac{n}{b(n)}(nT(n) - 1 - c(n)) \stackrel{d}{\Longrightarrow} \frac{1}{\mu^2}Y_2(1)$$

It remains to discuss the right hand side $II$ given by:

$$II = \frac{xnc(n)}{nc(n) - xb(n)}.$$

Clearly we have
$$II = \frac{x}{1 - xb(n)/(nc(n))}$$
and we have
$$\frac{b(n)}{nc(n)} = \frac{1}{c(n)} \frac{b(n)}{d(n)} \frac{d(n)}{n}$$
Since $X^2 \in D(Y_2(1))$, we automatically have $b(n)/d(n) \to 0$. If $\mu_2 < \infty$, it follows that $d(n)/n \to \mu_2$, and then $c(n) \to c$, a finite constant. In this case we conclude that $b(n) = o(1)nc(n)$. On the other hand, if $\mu_2 = \infty$, then $d(n)/n \to \infty$ and $nc(n)/d(n) \to 1$. Now it also follows that $b(n) = o(1)nc(n)$. We conclude that
$$I \to P(\frac{1}{\mu^2} Y_2(1) \leq x)$$

This proves the result.

(iv) To prove this result, as in the proof of (iii) we write
$$I = P(\frac{nc^2}{b(n)}(\frac{1}{c} - \frac{t^2(n)}{n}) \leq x)$$
and now we find that
$$I = P(\frac{n}{b(n)}(nT(n) - \frac{\mu_2}{\mu^2}) \leq \frac{ncx}{nc - xb(n)})$$

Using Theorem 4(iv) and $b(n)/n \to 0$, we obtain that
$$I \to P(\frac{1}{\mu^2} Y_2(\alpha/2) \leq x)$$

and the proof of the result.

(v) Similar as the proof of part (iii). ∎

# 5   Concluding remarks

1) There are many statistics that use higher sample moments. One can analyze domains of attraction of the random vector $(X, X^2, ..., X^k)$ and then apply the results to obtain weak limit theorems for these statistics.

2) The coefficient of variation and the sample dispersion are widely used measures of variation. For applications in the context of insurance and actuarial risk, we refer to [1], [13]. See also [17] and the references given there. In portfolio theory, the very popular ratio of Sharpe turns out to be given by $1/SV(n)$, cf. [23], [14]. The coefficient of variation is also used as a performance measure in queueing systems and in simulation, cf. [15].

# References

[1] Albrecher, H. and Teugels, J.L., 2004, Asymptotic analysis of measures of variation. Technical Report 2004-042, EURANDOM, T.U. Eindhoven, The Netherlands.

[2] Bentkus,V., Jing, B.Y., Shao, Q.M. and Zhou, W., 2006, Limiting distributions of the non-central $t$-statistic and their applications to the power of $t$-tests under non-normality. *Bernoulli* **13**:2, 346-364. MR2331255

[3] Bingham, N.H, Goldie, C.M. and Teugels, J.L., 1989, *Regular Variation*. Cambridge University Press. MR1015093

[4] Chistyakov, G.P. and Götze, F., 1996, Limit distributions of studentized means. *Ann. Probab.* 32, 28-77. MR2040775

[5] Feller, W. (1971). *An Introduction to Probability Theory and its Applications*, Vol. II (2nd edition). John Wiley and Sons, New York. MR0270403

[6] Geluk, J. and de Haan, L. (1987). *Regular Variation, Extensions and Tauberian Theorems*. CWI Tracts 40, Centre for Mathematics and Computer Science, Amsterdam, The Netherlands. MR0906871

[7] Giné, E., Götze, F. and Mason, D.M., 1997, When is the Student $t$-statistic asymptotically standard normal? *Ann. Probab.* **25**, 1514-1531. MR1457629

[8] Greenwood, P. and Resnick, S., 1979, A bivariate stable characterization and domains of attraction. *J. Multiv. An.* **9**, 206-221. MR0538402

[9] de Haan, L. and Omey, E., 1983, Integrals and derivatives of regularly varying functions in $R^d$ and domains of attraction of stable distributions II. *Stoch. Proc. and Their Appl.* **16**, 157-170. MR0724062

[10] de Haan, L, Omey, E. and Resnick, S., 1984, Domains of attraction and regular variation in $R^d$. *J. Muliv. An.* **14**, 17-33. MR0734097

[11] de Haan, L. and Ferreira, A., 2006, *Extreme Value Theory: An Introduction*. Springer-Verlag, New York. MR2234156

[12] Jessen, A.H. and Mikosch, T., 2006, Regularly varying functions. *Publ. Inst. Math. Béograd (N.S)*, **80 (94)**, 171-192. MR2281913

[13] Kaas,R., Goovaerts, M., Dhaene, J. and Denuit, M., 2001, *Modern Actuarial Risk Theory*. Kluwer Academic Publishers, Boston.

[14] Knight, J. and Stachell, S., 2005, A re-examination of Sharpe's ratio for log-normal prices. *Appl. Math. Finance* **12**:1, 87-100.

[15] Krishnamurthy, A. and Suri, R., 2006, Performance analysis of singlestage kaban controlled production systems using parametric decomposition. *Queueing Systems* **54**, 141-162. MR2268058

[16] Ladoucette, S.A. and Teugels, J.L., 2006, Asymptotic analysis for the ratio of the random sum of squares to the suare of the random sum with applications to risk measures. *Publ. Inst. Math., Nouv. Sér.* **80 (94)**, 219-240. MR2281916

[17] Ladoucette, S.A., 2007, Analysis of Heavy-Tailed Risks. Ph.D. Thesis, Catholic University of Leuven.

[18] Logan, B.F., Mallows, C.I., Rice, S.O. and Shepp, L.A., 1973, Limit distributions of self-normalized sums. *Ann. Probab.* **1**, 788 - 809. MR0362449

[19] Petrov, V.V., 1995, *Limit theorems of probability theory*. Oxford Science Publications. Oxford Studies in Probability Series. Clarendon Press - Oxford. MR1353441

[20] Resnick, S., 1986, Point processes, regular variation and weak convergence. *Adv. Appl. Prob.* **18**, 66-138. MR0827332

[21] Rvaceva, E., 1962, On the domains of attraction of multidimensional distributions. *Selected Transl. Math.Stat. Prob.* **2**, 183-207. MR0150795

[22] Seneta, E., 1976, *Functions of regular variation*. Lecture Notes in Mathematics 506. Springer, New York. MR0453936

[23] Sharpe, W.F., 1966, Mutual fund performance. *J. Business* **39**, 119-138.